\documentclass[11pt,reqno,twoside]{amsart}

\usepackage{bm}
\usepackage[english]{babel} 
\usepackage{esint}
\usepackage[T1]{fontenc} 
\usepackage{amsmath}
\usepackage{amssymb} 
\usepackage{float}
\usepackage{amsfonts}
\usepackage{mathtools}
\usepackage{xcolor}
\usepackage[utf8]{inputenc}
\usepackage[mathcal]{eucal} 
\usepackage{enumerate}
\usepackage{enumitem}
\usepackage{amsthm} 

\usepackage[totalwidth=14cm,bottom=1.25in, top=1.25in, hmarginratio=1:1]{geometry}

\newcommand{\R}{\mathbb{R}}

\newcommand{\oct}{\mathbb{O}'}
\newcommand{\K}{\mathcal{K}}

\newcommand{\Sp}{\mathrm{Sp}}

\newcommand{\SL}{\mathrm{SL}}
\newcommand{\Mix}{\mathrm{Mix}}

\newcommand{\Area}{\mathrm{Area}}

\newcommand{\hyp}{\mathrm{hyp}}

\newcommand{\zbar}{\overline{z}}

\newtheorem{defi}{Definition}[section] 
\newtheorem{thm}[defi]{Theorem}
 
\newtheorem{lemma}[defi]{Lemma}

\newtheorem{prop}[defi]{Proposition}

\newtheorem{bigthm}{Theorem}

\numberwithin{equation}{section}
\theoremstyle{definition}
\newtheorem{rmk}{Remark}[section]

\DeclareMathAlphabet{\mathpzc}{OT1}{pzc}{m}{it}

\usepackage[citecolor=teal,linkcolor=blue,urlcolor=cyan]{hyperref} 
\hypersetup{colorlinks}

\begin{document}

\title{Compactifying the $G_2'$-Hitchin component via almost-complex curves}
\author{Parker Evans, Charles Ouyang and Andrea Tamburelli}

\begin{abstract} 
By studying limits of induced metrics of equivariant almost-complex curves in the pseudosphere $\mathbb{S}^{2,4}$, we compactify the $G_{2}'$-Hitchin component.
\end{abstract}

\maketitle
\setcounter{tocdepth}{1}
\vspace{-1cm}

\section*{Introduction}

Higher rank Teichm{\"u}ller theory concerns representations of the fundamental group of a surface $S$ to Lie groups $G$ of real rank at least two. When $G$ is a semisimple Lie group, representations of $\pi_1S$ are encoded through the non-Abelian Hodge correspondence by the data of a Higgs bundle on a Riemann surface. At the heart of this correspondence is the existence and uniqueness of equivariant harmonic maps from the universal cover of the Riemann surface to the Riemannian symmetric space of $G$ \cite{Donaldson_selfdual, Corlette}, and the solution to a coupled system of elliptic PDE called \emph{Hitchin equations}, which govern harmonicity \cite{Hitchin_selfdual, Simpson_Higgs}. When the Lie group is real rank two, split, and simple, then for a unique choice of Riemann surface, the harmonic map is additionally conformal, thus giving a distinguished equivariant minimal surface in the symmetric space \cite{Labourie_cyclic}. In this case, the Hitchin equations simplify to either a single scalar equation or a coupled system of two equations.

In recent years, there has been tremendous interest in compactifications of distinguished components of the character variety, including the \emph{Hitchin components}, introduced by Hitchin in \cite{Hitchin_Teichmuller} using Higgs bundles. Most of these compactifications have tried to some extent to generalize the Thurston compactification of Teichm\"uller space, corresponding to  $G=\text{PSL}_{2} \mathbb{R}$. Of note, Parreau's compactification procedure \cite{parreau_compact} allows for any non-compact semisimple group $G$, and the limits obtained are actions on affine buildings, mirroring actions on $\mathbb{R}$-trees in the case of $\text{PSL}_{2} \mathbb{R}$ \cite{bestvina1988degenerations}. More recently Burger--Iozzi--Parreau--Pozzetti \cite{BIPPcompact, BIPP25} and Flamm \cite{Flam26} have constructed a real length spectrum compactification procedure, which uses the algebraic data of the representation. For compactifications using geodesic currents developed by Bonahon \cite{Bonahon_currents}, see \cite{Bonahon_currents, OT, OT_Sp4, Rei25}.

Even in rank two, there remains active work in understanding finer perspectives of limits, though at present many of these techniques only work for rays or sequences asymptotic to rays \cite{LTW_SL3, Rei26, WMa}. From a compactification perspective, this only furnishes a partial one, as often the case of degenerating conformal structures is difficult to study (see for example \cite{Loftin_compactify, Loftin_neckseparation}).

In this present case, we consider Hitchin representations in the exceptional split real Lie group $G_2'$. 
Such representations $\rho:\pi_1S\rightarrow G_2'$ 
admit equivariant almost-complex curves in the pseudosphere $\mathbb{S}^{2,4}$. In fact, the group $G_2'$ is realized as the automorphism group of the pseudo-Riemannian almost-complex manifold $\mathbb{S}^{2,4}$.  
Baraglia noted the existence of such equivariant curves $\nu: \tilde{\Sigma} \rightarrow \mathbb{S}^{2,4}$ for $G_2'$-Hitchin representations in \cite{Baraglia_thesis}. Due to work by Labourie in \cite{Labourie_cyclic} and Collier-Toulisse \cite{CT24}, each $G_2'$-Hitchin representation admits a canonical equivariant \emph{alternating} curve.

Baraglia's construction can be summarized as follows. 
From a $G'_{2}$-Hitchin representation $\rho$, via Labourie's preferred conformal structure $\Sigma$ on $S$, one can consider a Higgs bundle $(\mathcal{E},\Phi)$ of form \eqref{HitchinSection} in the $G_2'$-Hitchin section associated to a holomorphic sextic differential on $\Sigma$. The tautological section
$s \in \Omega^0(\Sigma, \mathcal{O})$ of the vector bundle $\mathcal{E}$, given by $s(p) =1 \in \mathcal{O}^{\R}_p$, corresponds, under the flat structure $\mathcal{E}^{\R} \cong \tilde{S}\times_{\rho} \R^{3,4}$, to a $\rho$-equivariant map $\nu: \tilde{\Sigma} \rightarrow \mathbb{S}^{2,4}$. This map $\nu$ is the almost-complex curve of interest. One can think of this curve $\nu$ as the $G_2'$-analogue of the equivariant affine sphere naturally associated to an $\mathrm{SL}(3,\R)$-Hitchin representation by Labourie-Loftin. Indeed, the latter was also  reinterpreted in similar fashion via Higgs bundles in Baraglia's thesis \cite{Baraglia_thesis}. 
In the $G_2'$-setting, the generalized Gauss map of the curve $\nu$ is precisely the unique $\rho$-equivariant conformal harmonic map to the $G'_{2}$-symmetric space \cite{Labourie_cyclic, CT24}. Hence, these almost-complex curves are the canonical geometric objects to study. 

We obtain a compactification of the $G'_{2}$-Hitchin component by studying the degeneration of the induced metrics of these curves, including the case that the conformal structure degenerates. 
Understanding the limits of these metrics along degenerating sequences in the character variety furnishes an explicit geometric approach to compactification procedures in higher Teichm{\"u}ller theory. However, it is important to note that this concrete approach comes at the expense of generality, as it works only for select Lie groups $G$, namely when there is a natural geometric object associated to a $G$-Hitchin representation.

We now describe the geometry of the degenerations of such curves. 
We show that properly rescaled, regions of the surface become flat and converge to a singular flat metric coming from a meromorphic sextic differential on a Riemann surface of lower topological complexity. On the complement, the rescaled induced metrics  converge to a measured lamination.
These hybrid objects are known as \emph{mixed structures}, introduced by Duchin-Leininger-Rafi in \cite{DLR_flat}. 
Henceforth, the space of projectivized mixed structures will be denoted $\mathbb{P}\text{Mix}_{6}(S)$. 
We use the language of geodesic currents developed by Bonahon \cite{Bonahon_currents} to give a convenient topology. 
In this way, we have a map $\textnormal{Hit}(S,G_{2}') \rightarrow \mathbb{P}\mathrm{Curr}(S)$ from the Hitchin component to the space of projective geodesic currents, which is compact, by reinterpreting the induced metric of the curve associated to a Hitchin representation. 
Our main result is a precise description of the boundary objects in the resulting compactification. 

\begin{bigthm} 
The boundary of $\textnormal{Hit}(S,G_{2}')$ in the space of projective geodesic currents is exactly $\mathbb{P}\textnormal{Mix}_{6}(S)$. 
\end{bigthm}
We note that when restricted to the copy of Teichm\"uller space given by the Fuchsian locus, which is embedded in $\mathbb{P}\mathrm{Curr}(S)$, this restricts exactly to Bonahon's realization of Thurston's compactification. 

In the case of the other split real simple Lie groups of rank two, previous work of the second and third named authors \cite{OT, OT_Sp4} have provided a similar compactifications using geometric objects related to $\SL_{3} \mathbb{R}$ and $\Sp_{4} \mathbb{R}$-Hitchin representations, respectively. 
Thus, this paper provides closure to this list.

\subsection*{Acknowledgments}
A.T. acknowledges support from the Italian Ministry of University and Research (MUR) through the FIS3 project “Geometry and Topology of higher Teichm\"uller spaces”, grant no. FIS-2024-05142. 

\subsection*{AI declaration}
No AI tools were used whatsoever in the production of this article.

\section{Background material}

\subsection{$G_2'$ and the Split Octonions}\label{subsec:G2}
The Lie group $G_2'$, the split real simple (adjoint) group of type $G_2$, can be described as follows.

We begin with the pseudo-Euclidean space $\R^{3,4} =(\R^7,Q_{3,4})$, where $Q=Q_{3,4}$ is a quadratic form of signature $(3,4)$. The space $\R^{3,4}$ admits a \emph{cross-product} $\times_{3,4}: \R^{3,4}\times \R^{3,4}\rightarrow \R^{3,4}$ that is bilinear, alternating, and normalized relative to $Q_{3,4}$ in similar fashion to the familiar cross product in Euclidean space $\mathbb{E}^3$, see \cite{DE26b}. 
The group $G_2'$ can be characterized as the automorphism group $G_2'= \mathrm{Aut}(\R^{3,4},\times_{3,4})$. Note further that
$\times_{3,4}$ is unique up to isomorphism, as from such a product one can recover the \emph{split octonions} $\oct$ \cite{Fon18}. The group $G_2'$ happens to preserve both the volume form on $\R^7$ and the quadratic form $Q$, making $G_2'$ a Lie subgroup of $SO(3,4)$. See \cite[Section 2]{DE26b} for further details and background.

The pseudosphere $\mathbb{S}^{2,4}$ can then be defined as the positive quadric in $\R^{3,4}$: 
\[ \mathbb{S}^{2,4} = \{ v \in \mathbb{R}^{3,4} : Q(v) =+1\}. \]
The cross product $\times_{3,4}$ induces an almost-complex structure $J$ on $S^{2,4}$, which can be seen as 
$J_u(v) := u\times_{3,4} v$, where the tangent space $T_u\mathbb{S}^{2,4}$ is identified with $u^{\bot_Q}$. This almost-complex structure $J$ is non-integrable \cite{CT24}. 
Recall that for an almost-complex manifold $(X, J)$ and a Riemann surface $(S, j)$, an immersion $f: (S,j) \to (X,J)$ is called an \emph{almost-complex curve}, or a $J$-holomorphic curve, if $df \circ j = J \circ df$. 

The pseudosphere carries a signature $(2,4)$ pseudo-Riemannian metric, inherited from $Q$ by restriction. Thus, almost-complex curves in $\mathbb{S}^{2,4}$ can be of different types, depending on the signature of their tangent and normal spaces. We shall consider \emph{alternating} curves, named by Nie in \cite{Nie24}, which have timelike tangent $T$, spacelike normal $N$, and timelike binormal spaces $B$. Equivalently, one can consider the anti-isometric partner $\mathbb{H}^{4,2}$ of $\mathbb{S}^{2,4}$ as in \cite{CT24} and \cite{Nie24}, which switches the signature of each such subspace. 

\subsection{The $G_2'$-Hitchin Component}\label{Subsec:G2Hitchin}
Let $S$ now be closed of genus $g \geq 2$. 
In \cite{Hitchin_Teichmuller}, Hitchin parametrized the $G_2'$-Hitchin component $\mathrm{Hit}(S,G_2')$ in the $G_2'$-character variety using the vector bundle $H^0(\Sigma_0, \K^2_{\Sigma_0})\oplus H^0(\Sigma_0,\K^{6}_{\Sigma_0})$ of holomorphic quadratic and sextic differentials on a fixed Riemann surface $\Sigma_0 = (S,j_0)$ on $S$. 

Hitchin's original work treated all split, real, simple groups $G$ of adjoint type. However, the $G_2'$-Hitchin case was not considered in particular until work by Baraglia in \cite{Baraglia_thesis}.  
There, Baraglia provided a construction for obtaining equivariant alternating almost-complex curves in $\mathbb{S}^{2,4}$ via Hitchin's parametrization. 
These curves were recently revisited and defined in greater generality in \cite{Evans24, Nie24, CT24, DE26b}. 

We now describe Hitchin's parametrization and Baraglia's construction in a bit more detail. 
Let $\Sigma$ be a compact Riemann surface of genus at least two.
We will partially describe the $G'_{2}$-Hitchin section, isolating the contribution from sextic differentials. Later, we will discuss how work of Labourie promotes these Higgs bundles to a full description of the Hitchin component by also changing the Riemann surface. 

Denote by $\mathcal{K}$ the holomorphic cotangent bundle of $\Sigma$. We consider the following Higgs bundles $(\mathcal{E},\Phi)$ on $\Sigma$:
\begin{align}\label{HitchinSection}
   \mathcal{E} = \bigoplus_{i=3}^{-3}\mathcal{K}^i, \qquad  \Phi = \begin{pmatrix} 0&  & & & & q_6&\\
				     \sqrt{3} & 0&  & & & & q_6\\
			  	      & \sqrt{5}&0 & & & &\\
			 	      & & -\sqrt{-6}&0 && &\\
				      & & & -\sqrt{-6}& 0& &\\
				      & & & &\sqrt{5} &0&\\
				      & & & & & \sqrt{3}&0\\
				     \end{pmatrix}  \ .				
	\end{align} 
Here, the Higgs field $\Phi = \Phi(q_6)$ is constructed from a holomorphic sextic differential $q_6 \in H^{0}(\Sigma, \mathcal{K}^{6})$ on $\Sigma$. The constants in the Higgs field $\Phi$ arise from $G_2'$ Lie theory; see \cite{Evans24}. We will write $q=q_6$ for simplicity going forward. For a discussion of the $G_2'$-structure on this Higgs bundle, see \cite{Eva25}, or \cite[Section 3.2]{DE26b}, which includes placing a non-degenerate holomorphic bilinear form $\mathcal{Q}$ and cross-product $\times_{\mathcal{E}}$ on $\mathcal{E}$, which $\Phi$ respects. 

By \cite{Baraglia_cyclic, CT24}, the solution to Hitchin's self-duality equations for these Higgs bundles $(\mathcal{E}, \Phi)$ is unique and given by a diagonal metric of the form
\begin{align}\label{HitchinMetric}
	h( \mathbf{x}, \mathbf{y} ) = \mathbf{\overline{x}}^T \, \mathsf{diag} (r^{-1}s^{-1}, r^{-1}, s^{-1}, 1, s, r, rs)\, \mathbf{y} =  \mathbf{\overline{x}}^T \,H \, \mathbf{y}.
\end{align} 
Here, $s$ is a metric on $\mathcal{K}^{-1}$ and $r$ a metric on $\mathcal{K}^{-2}$. The adjoint with respect to $h$ of the Higgs field $\Phi$ is given by
	\begin{align}
				  \Phi^{*h} &= H^{-1} \overline{\Phi}^T H = \begin{pmatrix} 0& \sqrt{3}\, s & & & & &\\
				      & 0& \sqrt{5} \frac{r}{s}& & & & \\
			  	      & &0 &\sqrt{-6}s  & & &\\
			 	      & & &0 &\sqrt{-6}s & &\\
				      & & & & 0&\sqrt{5} \frac{r}{s}&\\
				      \frac{ \overline{q}}{ r^2s }& & & & &0& \sqrt{3} \,s \\
				      & \frac{ \overline{q}}{ r^2s } & & & &&0\\
				     \end{pmatrix} 	.
	\end{align} 

In a local $z$-coordinate on $\Sigma$, Hitchin's equations for the flatness of the connection $\nabla = \nabla^{h}+\Phi+\Phi^{*h}$ can be written as the following system of PDE's for the flat Laplacian $\Delta = \partial_z \circ \partial_{\zbar}$:
 \begin{align}\label{HitchinsEquations_rs}
 	\begin{cases} 
 		\Delta \log r = 5 \frac{r}{s} - 3s- \, \frac{|q|^2}{r^2 s}, \, \\
		\Delta \log s = 6s-  5\frac{r}{s} \ .
	\end{cases} 
 \end{align} 
 Substituting $x_1= \log r + \frac{1}{2}\log s$, $x_2 = \frac{1}{2}\log s$, define $\psi = (\psi_1, \psi_2)$ by 
 \[ \begin{cases}
 	e^{\psi_1} \sigma^{5/2} = e^{x_1}, \\
	e^{\psi_2} \sigma^{1/2} = e^{x_2} .
 \end{cases} \] 
The functions $\psi_1, \psi_2: \Sigma \rightarrow \R$ are globally defined. 
 Hence, we may write the Hitchin equations in the following global form:
\begin{align}\label{G2Hitchin_GeneralMetric} 
	\begin{cases}  
		2 \Delta_\sigma \psi_1 = 5e^{\psi_1 - 3 \psi_2} -2 |q|^2_{\sigma} \, e^{-2\psi_1} + \frac{5}{2} \kappa_{\sigma}\ , \\
		2 \Delta_\sigma \psi_2 = -5e^{\psi_1 - 3\psi_2} + 6 e^{2\psi_2} + \frac{1}{2} \kappa_\sigma \,. \end{cases} 
\end{align} 
Here, $\sigma=\hat{\sigma}|dz|^2$ is a conformal metric on $\Sigma$ with curvature function $\kappa_{\sigma}= -2\frac{1}{\hat{\sigma}}(\log\hat{\sigma})_{z\zbar}$, and $|q|_{\sigma}^2=\frac{q \overline{q}}{\sigma^6}$. 
In what follows, we will set $\sigma=\sigma_{\mathrm{hyp}}$ to be the unique hyperbolic conformal metric on $\Sigma$.  

\begin{rmk}\label{remk:FuchsianHitchin}
Note that $\psi_1=0$, $\psi_2=0$ is \emph{not} the solution to \eqref{G2Hitchin_GeneralMetric} when $\sigma = \sigma_{\hyp}$ and $q=0$. Instead, one finds $e^{2\psi_2}\equiv \frac{1}{2}$ and $e^{\psi_1-3\psi_2}\equiv \frac{1}{2}$. 
\end{rmk}

Baraglia's construction of almost-complex curves from the Higgs bundles in the $G_2'$-Hitchin section was as follows. 

\begin{prop}[Associated almost-complex curve]\label{Prop:ACCurve}
Let $\rho:\pi_1S\rightarrow G_2'$ be a Hitchin representation
and $(\mathcal{E},\Phi)$ be the associated Higgs bundle of the form \eqref{HitchinSection} on a Riemann surface $\Sigma=(S,j)$. Then the tautological section $s \in \Omega^{0}(\Sigma, \mathcal{O})$ corresponds 
to a $\rho$-equivariant almost-complex curve $\nu:\tilde{\Sigma} \rightarrow \mathbb{S}^{2,4}$. 
\end{prop}

The rough idea is as follows. The section $s$ is real, namely valued in $\mathcal{E}^{\R}$, has norm $+1$ under the $\nabla$-parallel quadratic form $\mathcal{Q}$ of signature $(3,4)$ and thus $\sigma$ gives a $\rho$-equivariant map $\nu:\tilde{\Sigma} \rightarrow \mathbb{S}^{2,4}$. Then under the $G_2'$-structure on $\mathcal{E}^{\R}$, one finds $\nu$ is almost-complex and alternating. See \cite[Section 3.2.1]{DE26b} or \cite{CT24}.\medskip 

\begin{rmk}
The almost-complex curve $\nu$ is timelike in $\mathbb{S}^{2,4}$. We shall refer to 
the Riemannian metric $-\nu^*g_{\mathbb{S}^{2,4}}$ as the \emph{induced metric} of the almost-complex curve $\nu$. Up to a constant multiplicative factor, this metric is expressed as $g=e^{2\psi_{2}}\sigma$ in terms of $(\psi_1,\psi_2)$ from \eqref{G2Hitchin_GeneralMetric}.
\end{rmk}

In Hitchin's work \cite{Hitchin_Teichmuller}, the Riemann surface $\Sigma$ was fixed, and both $q_6\in H^0(\K^6_{\Sigma})$ and a holomorphic quadratic differential $q_2 \in H^0(\K^2_{\Sigma})$ varied in order parametrize the full $G'_{2}$-Hitchin component. We have suppressed the contribution from $q_2$ above in the Higgs field $\Phi$. 
However, Labourie's work in \cite{Labourie_cyclic} gives a parametrization of $\mathrm{Hit}(S,G_2')$ independent of a choice of basepoint in Teichm\"uller space. 
Indeed, combining Hitchin's parametrization of the $G_2'$-Hitchin component from \cite{Hitchin_Teichmuller} with Labourie's choice of canonical Riemann surface structure from \cite{Labourie_cyclic}, one obtains a diffeomorphic identification between $\mathrm{Hit}(S,G_2')$ and the total space of the bundle $\mathcal{Q}_6$ of holomorphic sextic differentials over Teichm\"uller space. 
Similarly, Collier-Toulisse gave in \cite{CT24} a homeomorphic interpretation of $\mathrm{Hit}(S,G_2')$ with a component of the moduli space of equivariant alternating almost-complex curves $(\nu,\rho)$, which eventually arrives back at the same description of Higgs bundles. 

\subsection{Geodesic Currents} 

The coarsest measure of a metric on a closed surface $S$ is via its marked length spectrum. If the metric $g$ has a length minimizing curve in each free homotopy class, one can consider its marked length spectrum given by the collection of pairs $\{ [\gamma], l_{g}([\gamma])\}_{\gamma \in \mathcal{C}(S)}$, where $l_{g}([\gamma])$ is the length, and $\mathcal{C}(S)$ is the set of free homotopy classes of closed loops in $S$. One can use geodesic currents (see \cite{Bonahon_currents}), which are $\pi_{1}(S)$-equivariant Radon measures on pairs of distinct unoriented bi-infinite geodesics on the universal cover of a hyperbolic surface $(S, \sigma)$. The space of geodesic currents Curr$(S)$ is defined via a weak-* topology, and the space of projectivized geodesic currents is its projectivization. In particular, $\mathbb{P}\mathrm{Curr}(S)$ is compact.

Nonpositively curved metrics and measured laminations may be realized as geodesic currents. On the space of geodesic currents, there is a continuous pairing $i(\cdot, \cdot)$, for which we isolate three important properties for this paper: 
\begin{enumerate}[label=(\roman*), itemsep=1ex]
    \item $i(L_{g}, L_{\gamma}) = l_{g}([\gamma])$, where $L_{g}$ is a geodesic current representing a metric $g$, and $L_{\gamma}$ is a closed curve.
    \item $i(L_{g}, L_{g}) = \frac{\pi}{2} \text{Area}(S,g) $.
    \item $i(\eta, \eta) = 0$ if and only if $\eta$ represents a measured lamination. 
\end{enumerate}

A particular important subclass of nonpositively curved metrics are obtained by holomorphic $n$-differentials on a compact Riemann surface $\Sigma$. If $q$ denotes a holomorphic $n$-differential on $\Sigma$, then $|q|^{2/n}$ is a singular flat metric. Away from the zeros of $q$, holomorphicity implies the metric is flat; at an order $k$ zero of $q$, there is a conical singularity of angle $2\pi + \frac{2 \pi}{k}$. Another important class of geodesic currents are the \emph{$n$-mixed structures}. For a surface $S$, one takes a collection of disjoint incompressible subsurfaces $S_{\alpha}$, and chooses a finite area meromorphic $n$ differential $q$, so that the boundary curves all have $|q|^{2/n}$-length zero. On the complement of all such subsurfaces, choose a measured lamination $\lambda$. Then an $n$-mixed structure is a geodesic current $\eta = L_{|q|^{2/n}} + \lambda$. 
Note that a mixed structure allows for the two extreme situations of a holomorphic differential $q$ on all of $\Sigma= (S,j)$ and of a pure measured lamination. 
Zhang has shown in her thesis \cite{ZZhang} that limits of degenerating sequences of holomorphic $n$-differential metrics converge projectively to $n$-mixed structures. The case of $n=2$ was first shown in \cite{DLR_flat}, and $n=3,4$ in \cite{OT}, \cite{OT_Sp4}.
We shall denote $\mathbb{P}\Mix_6(S)$ for the space of 6-mixed structures, viewed in the space $\mathbb{P}\mathrm{Curr}(S)$ of projective geodesic currents. 

\section{Main results}

The proof of the main theorem follows the strategy developed in \cite{OT,OT_Sp4}. Given a Hitchin representation $\rho$, we compare the induced metric $g$ of the almost-complex curve $\nu$ with the singular flat metric $|q|^{1/3}$ associated to the sextic differential. The maximum principle applied to Hitchin’s equations, shows that $g$ has negative curvature (Proposition \ref{prop:curvature}), and lower bounds for $g$ in terms of both the conformal hyperbolic metric (Proposition \ref{prop:Domination}) and \(|q|^{1/3}\) (Lemma \ref{lower}). Integrating these equations shows that the area of $g$ asymptotically behaves like a fixed multiple of the flat area (Proposition \ref{prop:AreaBound}). These estimates allow us to analyze the behavior of the length spectrum of $g$ along divergent sequences of Hitchin representations. The key analytic ingredient is an exponential decay estimate for the difference between the conformal factors of the induced and flat metrics away from the zeros of the differential (Lemma \ref{exp}).\medskip

Throughout this section, we fix our notation as in the preliminaries. Given a Hitchin representation $\rho:\pi_1S\rightarrow G_2'$, we have a canonically associated Higgs bundle $(\mathcal{E},\Phi, \Sigma)$ as in \eqref{HitchinSection} from which we extract a $\rho$-equivariant curve $\nu:\tilde{\Sigma}\rightarrow \mathbb{S}^{2,4}$ via Proposition \ref{Prop:ACCurve}, which has induced Riemannian metric denoted by $g$. We shall write $\sigma=\sigma_{\hyp}$ for the unique conformal hyperbolic metric on the Riemann surface $\Sigma$. 

The following was proven in \cite[Theorem 8.4]{DaiLi24}. We include a short proof in our notation. As a consequence, the induced metrics may be regarded as geodesic currents. 

\begin{prop}\label{prop:curvature}
Let $\rho:\pi_1S\rightarrow G_2'$ be Hitchin and $g$ be the induced metric of the associated $\rho$-equivariant curve $\nu:\tilde{\Sigma}\rightarrow \mathbb{S}^{2,4}$. Then $g$ is negatively curved. 
\end{prop}

\begin{proof}
We compute 
$$ \kappa_{g} = \frac{-2}{e^{2\psi_2}} \big ( \Delta_{\sigma}(2\psi_2) + \Delta_{\sigma} \log \sigma \big) = -12+10e^{\psi_1 -3\psi_2}.$$
Thus, $\kappa_{g} \leq 0 $ if and only if 
    \begin{align}\label{UpperBound1}
        e^{ \psi_1 - 5\psi_2} \leq \frac{6}{5}.
    \end{align} 
Consider $\alpha = \log( |q|^2_{\sigma} e^{-2\psi_1 - 2\psi_2} )$. 
Then one computes 
$$ \Delta_{\sigma} \alpha = 2|q|^2_{\sigma} e^{-2\psi_1}-6e^{2\psi_2} .$$
Hence, at a maximum $p$, 
$$ e^{\alpha(p)} = 2|q|^2_\sigma e^{-2\psi_1 -2\psi_2}(p) -6 = \frac{1}{e^{2\psi_2}}\Delta_{\sigma} \alpha(p) \leq 0 .$$
Thus, 
\begin{align}\label{UpperBound2}
    2|q|^2_{\sigma}e^{-2\psi_1 -2 \psi_2} \leq 6.
\end{align} 
Now, we prove \eqref{UpperBound1}. 
One directly computes 
\begin{align} \label{1minus5}
    2 \Delta_{\sigma}(\psi_1 -5\psi_2) = 30e^{\psi_1 - 3\psi_2}- 2|q|^2_{\sigma} e^{-2\psi_1}-  30 e^{2\psi_2} .
\end{align}

Applying the maximum principle, along with the inequality 
\eqref{UpperBound2}, 
$$30 e^{\psi_1 -5\psi_2} - 36 \leq 30 e^{\psi_1 -5\psi_2} - 30 - 2 |q|^2_{\sigma} e^{-2\psi_1 - 2\psi_2}= \frac{1}{e^{2\psi_2}} 2\Delta_{\sigma}(\psi_1-5\psi_2) \leq 0. $$
It follows that \eqref{UpperBound1} holds.

In fact, the strong maximum principle asserts that $\psi_1-5\psi_2 \equiv 6$ or $\psi_1-5\psi_2 < 6$ globally. In the former case, one finds that $q$ must be non-vanishing, which is impossible on a closed surface $\Sigma$. Hence, the inequality \eqref{UpperBound1} is strict, 
which completes the proof.
\end{proof} 

Next, we note that the space of induced metrics of such almost-complex curves realizes a circle quotient of the $G_2'$-Hitchin component. 
\begin{prop}
Let $\rho,\rho':\pi_1S\rightarrow G_2'$ be Hitchin, with associated sextic differentials $q, q'$, respectively. 
Then the induced metrics $g$ and $g'$ coincide if and only if $q'=e^{i \theta} q$. 
\end{prop}
\begin{proof}
If $q'=e^{i \theta} q$, then the statement in one direction follows immediately, as Hitchin's equations \eqref{G2Hitchin_GeneralMetric}, and hence the solutions $(\psi_1,\psi_{2})$ and $(\psi_1',\psi_{2}')$, are the same. Conversely, if the metrics coincide, then $\psi_{2} = \psi_{2}'$. Subtracting the second equations in \eqref{G2Hitchin_GeneralMetric}, one sees $\psi_{1} = \psi_{1}'$. Then $|q| = |q'|$ by the first equation in \eqref{G2Hitchin_GeneralMetric}. 
\end{proof}

We now define constants $c,d > 0$ uniquely solving the following system:
    \begin{align}\label{cdSystem}
 	\begin{cases} 5c^{5/6}d^{-1/2} - 2c^{-5/3} = 0 \\
	-5c^{5/6}d^{-1/2} +6d^{1/3} = 0. \end{cases} 
\end{align}

The exact solution is $d =\frac{5}{6\sqrt{3}}$ and $c = \left(\frac{2}{5}\right)^{2/5} d^{1/5}$.
We can use $c,d$ to define sub-solutions to \eqref{G2Hitchin_GeneralMetric}, similar to \cite{Evans24}.
We rewrite the system \eqref{G2Hitchin_GeneralMetric} in the form
$F(\bm{u}) = (F_1(\bm{u}), F_2(\bm{u}) )$ going forward. 

\begin{prop}\label{prop:Domination}
Let $\rho:\pi_1S \rightarrow G_2'$ be a Hitchin representation and $g=e^{2\psi_2}\sigma$ the rescaled metric of the associated $\rho$-equivariant curve. Then $g \geq \frac{1}{2}\sigma$.
\end{prop}

\begin{proof}
Let us set $\alpha = -2\psi_2$. Then observe that for $g= e^{2\psi_2}\sigma$, we have 
\[ \Delta_{g} \alpha = -\frac{1}{e^{2\psi_2}}\Delta_{\sigma}(2\psi_2)= 5e^{\psi_1-5\psi_2}-6+\frac{1}{2}e^{\alpha}\]
A similar argument to that of inequality \eqref{UpperBound1} shows that 
$e^{\psi_1-5\psi_2} \geq 1$ holds globally. 
At a global maximum $p$ of $\alpha$, we find  
\[ 0 \geq \Delta_g \alpha(p) \geq \frac{1}{2}e^{\alpha(p)}-1 .\] 
Thus, $1\geq \frac{1}{2}e^{\alpha}$, or equivalently, $e^{2\psi_2} \geq \frac{1}{2}$, holds globally.
\end{proof}

\begin{rmk}
In fact, the inequality in Proposition \ref{prop:Domination} is sharp. In the Fuchsian-Hitchin case, $g= \frac{1}{2}\sigma$. See also the related Remark \ref{remk:FuchsianHitchin}. 
\end{rmk}

The following sub-solution is similar to \cite[Lemma 5.2]{Evans24}. 

\begin{lemma}\label{subsolution}Let \(\Sigma=(S,J)\) be a Riemann surface, let \(q\in H^0(\Sigma,K_\Sigma^6)\), and let \(\sigma\) be the conformal hyperbolic metric on $\Sigma$. Then the functions
\begin{align}
        v_1 &= \begin{cases}
                \frac{5}{6} \log(c |q|_\sigma) & |q|_{\sigma} \geq \frac{1}{8c}, \\
                \frac{5}{2}\log(\frac{1}{2}) & \text{otherwise}. \\
                \end{cases} \\
         v_2 &= \begin{cases}
         \frac{1}{6} \log(d |q|_\sigma) & |q|_{\sigma} \geq             \frac{1}{8d} ,\\
                \frac{1}{2}\log( \frac{1}{2}) & \text{otherwise}. 
    \end{cases}
\end{align}
defines a locally Lipschitz weak subsolution $v= (v_1,v_2)$ of \eqref{G2Hitchin_GeneralMetric}.
\end{lemma}
\begin{proof}
Write equation \eqref{G2Hitchin_GeneralMetric} as \(\Delta_\sigma\psi_i=F_i(\psi_1,\psi_2)\), where
\[
\begin{aligned}
F_1(\psi_1,\psi_2)&=\frac{5}{2}e^{\psi_1-3\psi_2}-|q|_{\sigma}^2e^{-2\psi_1}-\frac{5}{4},\\
F_2(\psi_1,\psi_2)&=-\frac{5}{2}e^{\psi_1-3\psi_2}+3e^{2\psi_2}-\frac{1}{4}.
\end{aligned}
\]
Consider the constant vector
\[
u=\left(\frac52\log\frac12,\,
         \frac12\log\frac12\right).
\]
Since \(e^{u_1-3u_2}=e^{2u_2}=1/2\), we have
\[
F_1(u)=-32|q|_{\sigma}^2\le0=\Delta_\sigma u_1,
\qquad
F_2(u)=0=\Delta_\sigma u_2.
\]
Thus \(u\) is a subsolution.
Away from the zeros of $q$, define
\[
w=\left(\frac{5}{6}\log(c|q|_{\sigma}),\,
         \frac{1}{6}\log(d|q|_{\sigma})\right).
\]
Holomorphicity of \(q\) and hyperbolicity of $\sigma$ give
\[
\Delta_\sigma\log |q|_{\sigma}=-\frac32.
\]
Using the relations between $c$ and $d$ from \eqref{cdSystem}, we therefore obtain
 \begin{align}
     \Delta_{\sigma} w_1 &=-\frac{5}{4} = -\frac{5}{4}+ |q|_{\sigma}^{1/3}\left(\frac{5}{2}c^{5/6}d^{-1/2}-c^{-5/3}\right)= F_1(w_1, w_2), \\
     \Delta_{\sigma} w_2 &= -\frac{1}{4} = -\frac{1}{4}+ |q|_{\sigma}^{1/3}\left(-\frac{5}{2}c^{5/6}d^{-1/2}+3d^{1/3} \right)= F_2(w_1, w_2).
     \end{align}
Hence \(w\) is an exact solution away from the zeros of $q$.
Finally, the signs of the derivatives
\[
\frac{\partial F_1}{\partial \psi_2}
=-\frac{15}{2}e^{\psi_1-3\psi_2}<0,
\qquad
\frac{\partial F_2}{\partial \psi_1}
=-\frac52e^{\psi_1-3\psi_2}<0.
\]
imply that the componentwise maximum \(v=\max\{u,w\}\) is a weak subsolution (see for example \cite[Theorem 2.5]{TW24}.)
\end{proof}

We now show the metric on the almost-complex curve dominates (a constant multiple of) the 
flat metric $|q|^{1/3} =|q|_{\sigma}^{1/3} \sigma$. 
\begin{lemma}\label{lower}
The metric $g =e^{2\psi_2}\sigma$ satisfies 
\[ g \geq d^{1/3} |q|^{1/3},\]
where $d$ is from \eqref{cdSystem}. 
\end{lemma}

\begin{proof}
We shall consider the equations \eqref{G2Hitchin_GeneralMetric} with respect to the flat metric $|q|^{1/3}$. To this end, we first let $(x_1,x_2)$ and $(\psi_1,\psi_2)$ be given as in \S\ref{Subsec:G2Hitchin}. 
Define the pair $(\alpha_1, \alpha_2)$ so that 
$$\begin{cases}
    e^{x_1} = e^{\alpha_1}|q|^{5/6} \\ 
    e^{x_2} = e^{\alpha_2}|q|^{1/6}.
\end{cases}$$ 
In this notation, $g= e^{2u_2}|dz|^2$ locally. Following the reasoning in \S\ref{Subsec:G2Hitchin}, the functions $\alpha_1,\alpha_2$ are globally defined on $S \setminus Z(q)$, where $Z(q)$ is the zero set of $q$.

Then one finds that 
$$\begin{cases}
    e^{\psi_1} = e^{\alpha_1}|q|_{\sigma}^{5/6} \\ 
    e^{\psi_2} = e^{\alpha_2}|q|_{\sigma}^{1/6}.
\end{cases}$$ 
Now, the inequality $e^{\psi_1- 5\psi_2} \leq \frac{6}{5}$, verified in the proof of Proposition \ref{prop:curvature}, is equivalent to $e^{\alpha_1- 5\alpha_2} \leq \frac{6}{5}.$ Similarly, the inequality \eqref{UpperBound2} is equivalent to 
$$\alpha_1 + \alpha_2 \geq \frac{1}{2} \log\left( \frac{1}{3} \right)$$
Combining these two bounds, one finds 
$$ 6\alpha_2 \geq \log \left ( \frac{5}{6\sqrt{3}} \right) = \log d .$$
This inequality, which holds on $S \setminus Z(q)$, implies desired global inequality on $S$. 
\end{proof}

By integrating estimates for $g$, we obtain the following area comparisons.

\begin{prop}[Area Bounds]\label{prop:AreaBound}
Let $d = \frac{5}{6 \sqrt{3}}$ be as before, then
$$d^{1/3}||q|| \leq \mathrm{Area}(S, g) \leq d^{1/3}||q|| + \frac{7 \pi}{6} |\chi(S)|.$$
\end{prop}

\begin{proof}
    The lower bound follows from Lemma \ref{lower}. For the upper bound, we first claim that the following inequality holds globally:
\begin{align}\label{BC}
    |q|_{\sigma}^{2}e^{-2\psi_{1}-2\psi_2} \leq 3. 
\end{align}
If $q$ identically vanishes, the inequality vacuously holds. Thus, we suppose $q$ is nonzero. 
To verify the claim, apply Hitchin equations \eqref{G2Hitchin_GeneralMetric}. First, we compute:
\begin{align}\label{eq:Psi1Psi2}
    2\Delta_{\sigma} (\psi_{1} + \psi_{2} )= 6e^{2\psi_{2}} -2|q|_{\sigma}^{2}e^{-2\psi_{1}} + 3\kappa_{\sigma}.
\end{align}
Note that $\Delta_{\sigma} \log |q|_{\sigma}^{2} = 3 \kappa_{\sigma}$.
Now, set $\eta: =  \log\big(|q|_{\sigma}^{2}e^{-2(\psi_{1}+ \psi_{2})})$, 
which is a well-defined function away from the isolated and finite zero set $Z(q)$. 
We compute that 
\begin{align}\label{eq:Eta}
    \Delta_{g} \eta= 2e^{\eta} -6. 
\end{align}
Since $\eta \rightarrow -\infty$ near the zeros of $q$, the function $\eta$ attains a global maximum at some point $p$ on $S \setminus Z(q)$. 
Then \eqref{eq:Eta} implies $e^{\eta(p)}\leq 3$. Hence, $e^{\eta} \leq 3$ globally, proving the estimate \eqref{BC}. \medskip

Now, let us write 
\[d^{1/3} \frac{|q|_{\sigma}^{1/3}}{e^{2\psi_{2}}} = \left(\frac{5}{6} \cdot\frac{e^{\psi_{1}-3\psi_{2}}}{e^{2\psi_{2}}}\right)^{1/3} \left(\frac{1}{3} \cdot \frac{|q|_{\sigma}^{2}e^{-2\psi_{1}}}{e^{2\psi_{2}}} \right)^{1/6}.\]
From \eqref{UpperBound1} and \eqref{BC}, we have $0 < \left(\frac{5}{6} \cdot\frac{e^{\psi_{1}-3\psi_{2}}}{e^{2\psi_{2}}} \right)^{1/3} < 1$ and $0<\left( \frac{1}{3} \frac{|q|^{2}e^{-2\psi_{1}}}{e^{2\psi_{2}}} \right)^{1/6} <1$ respectively, so that
\[1- d^{1/3} \frac{|q|_{\sigma}^{1/3}}{e^{2\psi_{2}}} \leq \left(1- \left(\frac{5}{6} \cdot\frac{e^{\psi_{1}-3\psi_{2}}}{e^{2\psi_{2}}} \right)^{1/3} \right)  + \left(1 -\left( \frac{1}{3} \frac{|q|^{2}_{\sigma}e^{-2\psi_{1}}}{e^{2\psi_{2}}} \right)^{1/6} \right)\ .\]
But as $t^{\alpha}  \geq t$ for $0 \leq t \leq 1$ and $0 < \alpha < 1$, we have
$$ 1- d^{1/3} \frac{|q|_{\sigma}^{1/3}}{e^{2\psi_{2}}} \leq 1 - \frac{5}{6} \cdot\frac{e^{\psi_{1}-3\psi_{2}}}{e^{2\psi_{2}}} + 1 - \frac{1}{3} \frac{|q|_{\sigma}^{2}e^{-2\psi_{1}}}{e^{2\psi_{2}}}.$$
Hence, rewriting, we obtain
\begin{align}\label{upoint2}
e^{2 \psi_{2}} \leq d^{1/3}|q|_{\sigma}^{1/3} + \frac{5}{6}\bigg( \frac{6}{5} e^{2\psi_{2}} - e^{\psi_{1}-3\psi_{2}}\bigg)  +\frac{1}{3}(3e^{2\psi_{2}} - |q|_{\sigma}^{2}e^{-2\psi_{1}})\ .
\end{align}
Integrating the second equation of \eqref{G2Hitchin_GeneralMetric} and applying Gauss-Bonnet gives
\begin{align}\label{constant1}
    \int_{S} (6e^{2\psi_{2}} - 5 e^{\psi_{1}-3\psi_{2}}) dA_{\sigma} = \int_{S} -  \frac{1}{2} \kappa_\sigma \, dA_{\sigma} = \pi |\chi(S)|.
\end{align} 
Integrating the equation \eqref{eq:Psi1Psi2} yields
$$0 = 6\int_{S} e^{2\psi} dA_{\sigma} -2 \int_{S} |q|_{\sigma}^{2} e^{-2\psi_{1}} dA_{\sigma}+ 3\int_{S} \kappa_{\sigma}dA_{\sigma}\ .$$
Whence, 
\begin{align} \label{constant2}
    \int_{S} \big(3e^{2\psi_{2}} - |q|_{\sigma}^{2} e^{-2\psi_{1}}\big)dA_{\sigma} = 3 \pi|\chi(S)| 
\end{align}
so that integrating \eqref{upoint2} together with \eqref{constant1} and \eqref{constant2} gives the desired upper bound. 
\end{proof}

We require one last technical result, which provides asymptotics for the difference between the conformal factors of the induced metric and the associated flat metric.

\begin{lemma}\label{exp}
    If $q$ is a holomorphic sextic differential of sufficiently large norm, then the unit-area singular metric $|q_{0}|^{1/3} = |q|^{1/3}/||q||$ satisfies the following:
    if the $|q_{0}|^{1/3}$-ball of radius $r_0>0$ centered at a point $p$ is embedded and does not contain any zeros of $q$, then 
     \[ |\psi_{2}(p) - \frac{1}{6}\log (\, d|q(p)|_{\sigma})| \leq C e^{-c ||q||^{1/2} r_0},\]
     where $d = \frac{5}{6 \sqrt{3}}$, and the constants $c, C>0$ depend only on $r_0$. 
\end{lemma}

\begin{proof}
Define $u_{1}: = \psi_{1} - \frac{5}{6} \log(c |q|_{\sigma})$ and $u_{2}: = \psi_{2} - \frac{1}{6} \log(d|q|_{\sigma})$, which denotes the difference between the solution to the Hitchin equations \eqref{G2Hitchin_GeneralMetric} and the subsolutions given in the proof of Lemma \ref{subsolution}. In particular, we have $u_1, u_2 \geq 0$. Then
\begin{align*}
    \Delta_{\sigma} (u_{1} + u_{2}) &= \Delta_{\sigma} \psi_{1} + \Delta_{\sigma} \psi_{2} - \Delta_{\sigma} \log|q|_{\sigma} \\
    &= e^{2\psi_{2}} -|q|_{\sigma}^{2} e^{-2\psi_{1}} \\
    & = |q|_{\sigma}^{1/3} (3d^{1/3} e^{2u_{2}} -c^{-5/3}e^{-2u_{1}})\\
    &=3d^{1/3} |q|_{\sigma}^{1/3}(e^{2u_{2}}- e^{-2u_{1}})\ .
\end{align*}
Now, let us fix a point $p \in S \setminus Z(q)$, and a natural coordinate $z$ such that $q= dz^6$ on a neighborhood of $p$ that does not intersect $Z(q)$. Then we have
\begin{align}\label{deldelbar}
\partial_{z}{\partial}_{\zbar} (u_{1} + u_{2})= \Delta_{|q|^{1/3}} (u_{1} +u_{2}) &= 3d^{1/3} e^{-2u_{1}}(e^{2(u_{1}+u_{2})}-1)
\end{align}
As $u_{2}$ is the difference between the solution to the Hitchin equations and a subsolution, it is a subharmonic function. Using this together with Jensen's inequality on an embedded ball $B$ of $|q|^{1/3}$-radius $r$ gives
\begin{align*}
e^{2u_{2}} &\leq e^{\fint_{B} 2u_{2} \, dA_{q}} \leq \fint_{B} e^{2u_{2}} \, d A_{q}
{=} \frac{1}{\pi r^{2}}d^{-1/3} \int_{B} d A_{g} \leq d^{-1/3}\frac{\text{Area}(S, g)}{ \pi r^{2}}.
\end{align*}
The radii $r_0$ and $r$ of $B$, computed using the metrics $|q_{0}|^{1/3}$ and $|q|^{1/3}$, respectively, are related by $r= r_{0} ||q||^{1/2}$, so that we have the inequality
\begin{align}\label{Jensen}
    e^{2u_{2}} \leq d^{-1/3}\frac{\text{Area}(S,g)}{ \pi r_{0}^{2}||q||}.
\end{align}
Note that 
$$ u_{1} + u_{2} = u_{1} -5u_{2} + 6u_{2} = \psi_{1} - 5\psi_{2} - \frac{5}{6} \log\big( \frac{c}{d}|q|_{\sigma}\big) + 6u_{2}.$$
From the uniform upper bound \eqref{UpperBound1} and \eqref{Jensen} together with Proposition \ref{prop:AreaBound}, we obtain a uniform upper bound for $u_{1}+u_{2}$ on $B$ and consequently for $u_{1}$ as well.
Hence from \eqref{deldelbar} we now have,
$$ \partial_{z} {\partial}_{\zbar} (u_{1} + u_{2}) \geq a (e^{2(u_{1}+u_{2})}-1) \geq 2a  (u_{1} + u_{2}),$$
where $a$ is a fixed constant.
The solution to the Dirichlet boundary problem
$$ \partial_{z}{\partial}_{\zbar} \, \eta = a \eta ,\qquad \eta|_{\partial B} = A, $$
for constants $a >0$ and $A >0$, has a well-known (\cite[Eqs.~(10.25.1)--(10.25.2)]{DLMF}) radial solution given by a modified Bessel function of the first kind:
$$\eta(z) = A \frac{I_{0}(2\sqrt{2a} \,|z|)}{I_{0}(2 \sqrt{2a} r)} = A\frac{I_{0}(2\sqrt{2a} \,|z|)}{I_{0}(2 \sqrt{2a} ||q||^{1/2}r_{0})}, $$
where $I_0$ satisfies the asymptotic $I_{0}(x) \sim \frac{e^{x}}{2\pi x^{1/2}}$ as $x \to \infty$.
Hence, if $||q||$ is sufficiently large, 
$$ u_{2}(p)\leq u_{1}(p)+u_{2}(p) \leq Ce^{-c ||q||^{1/2} r_{0}},$$
as desired, for some positive constants $c,C$, which depend only on $r_0$.
\end{proof}

We now prove the main result.

\begin{thm}
Let $\{\rho_{n}\}$ be a sequence of $G'_{2}$-Hitchin representations leaving all compact sets. Denote by $(\mathcal{E}, \Phi_{n}, \Sigma_n)$ the corresponding $G'_{2}$-Higgs bundles and $g_{n}$ the induced metric of the resulting almost-complex curve. If $L_{g_{n}}$ is the associated geodesic current of $g_n$, then there are constants $t_{n} \to 0$ such that $t_{n}L_{g_{n}}$ subconverges to a mixed structure $\eta$. Conversely, given any mixed structure $\eta$, there is a renormalized sequence $t_{n}L_{g_{n}}$ arising from $G'_{2}$-Hitchin representations such that $t_{n}L_{g_{n}} \to \eta$. 

That is, the boundary $\mathrm{Hit}(S,G_2')$ in  $\mathbb{P}\mathrm{Curr}(S)$ is exactly $\mathbb{P}\Mix_6(S)$. 
\end{thm}

\begin{proof}
Let a sequence $\{\rho_{n} \}$ be given as in the hypotheses. 
As the space of projective geodesic currents is compact, we may suppose that
a renormalized sequence $t_{n}L_{g_{n}}$ subconverges to a geodesic current $\mu$. 
We now consider cases. 

\textbf{Case 1}: $||q_{n}||$ is uniformly bounded.
In this case, the underlying sequence of hyperbolic metrics $\{\sigma_{n}\}$ must diverge, so that passing to a subsequence there is a closed curve $\gamma$ on the surface whose $\sigma_{n}$-length diverges. 
By Proposition \ref{prop:Domination}, then the $g_{n}$-length of $\gamma$ also diverges. Hence, for the sequence $t_{n}L_{g_{n}}$ to converge, the constants $t_{n}$ must tend to zero. But Proposition \ref{prop:AreaBound} implies the area of $g_{n}$ is uniformly bounded, meaning $\Area(S,t_{n} g_{{n}}) \rightarrow 0$. Hence, the self-intersection $i(\mu,\mu)$ of $\mu$ must vanish, whence $\mu$ is a projectivized measured lamination.

\textbf{Case 2}: The sequence $||q_{n}||$ is unbounded. Hence, passing to a subsequence as needed, we can assume $\lim_{n \to \infty} ||q_{n}|| = \infty$. Writing $|q_{n}|^{1/3} = ||q_{n}||\cdot |q_{0,n}|^{1/3}$, where $|q_{0,n}|^{1/3}$ is a unit-area singular flat metric, define the sequence $\hat{L}_{g_{n}}:= \frac{1}{\sqrt{d^{1/3}||q_{n}||}} L_{g_{n}}$. If the length spectrum of this sequence is not uniformly bounded, then there are constants $t_{n} \to 0$ for which $t_{n}\hat{L}_{g_{n}}$ converges to $\mu$ and once again from Proposition \ref{prop:AreaBound}, $i(\mu, \mu) =0$ so that $\mu$ is a projectivized measured lamination. Hence, we may assume the length spectrum of $\hat{L}_{g_{n}}$ converges, and consequently from Lemma \ref{lower} the length spectrum of $|q_{0,n}|^{1/3}$ converges, say to $\eta$. But from \cite[Theorem 3.1]{ZZhang}, then $\eta$ is a mixed structure. 
We note that $\eta$ cannot be a measured lamination because the self-intersection of $\hat{L}_{g_{n}}$ is bounded away from $0$ by Proposition \ref{prop:AreaBound}.
Hence, $\eta$ is a mixed structure with nontrivial flat part. 
This means there is a decomposition of the surface $S$ given by a finite collection of disjoint simple closed curves, giving subsurfaces $S_{\alpha}$ with meromorphic sextic differential metrics $|q_{\infty}|^{\frac{1}{3}}$ with total area 1.

By Lemma \ref{exp}, away from small neighborhoods of zeros of $q_{\infty}$, the difference in conformal factors converges locally uniformly to 0, so that the length spectrum $L_{g_n}$ of the induced metrics $g_n$ on the subsurfaces $S_{\alpha}$ converge to that of $|q_{\infty}|^{1/3}$ (see Proposition 5.3 \cite{Charles_dPSL}). By Proposition \ref{prop:AreaBound}, the remaining area is zero.
Hence, $\eta= L_{|q_{\infty}|^{1/3}}+\lambda_{\infty}$, where $\lambda_{\infty}$ is a measured lamination, completing the proof in the forward direction.

The converse direction follows by first observing that the closure of flat metrics coming from holomorphic sextic differentials are precisely the mixed structures \cite{ZZhang}. As any flat metric $q$ may be obtained by looking at the sequence of induced metrics coming from the ray $tq$, the result follows.
\end{proof}

\bibliographystyle{alpha}
\bibliography{bs-bibliography-2-2}

\bigskip
\noindent \footnotesize \textsc{Department of Mathematics, Washington University, St. Louis}\\
\emph{E-mail address:} \verb|evans.p@wustl.edu|

\bigskip
\noindent \footnotesize \textsc{Department of Mathematics, Washington University, St. Louis}\\
\emph{E-mail address:}  \verb|ouyang@math.wustl.edu|

\bigskip
\noindent \footnotesize \textsc{Department of Mathematics, University of Pisa}\\
\emph{E-mail address:} \verb|andrea_tamburelli@libero.it|

\end{document}